\documentclass[12pt, reqno]{amsart}
\usepackage{amsmath, amsthm, amscd, amsfonts, amssymb, graphicx, color}
\usepackage[dvipsnames]{xcolor}
\usepackage[bookmarksnumbered, colorlinks, plainpages]{hyperref}
\hypersetup{colorlinks=true,linkcolor=red, anchorcolor=green, citecolor=cyan, urlcolor=red, filecolor=magenta, pdftoolbar=true}

\usepackage{tensor}

\newtheorem{defn}{Definition}[section]

\newtheorem{thm}{Theorem}[section]

\newtheorem{prop}[thm]{Proposition}
\newtheorem{cor}[thm]{Corollary}

\newtheorem{remark}[thm]{Remark}
\newtheorem{example}[thm]{Example}

\usepackage{enumitem}

\usepackage{fullpage}

\begin{document}
	
	\setcounter{page}{1}
	
	\title[FPT on BMS]{On Fixed Point Theorems in Bipolar Metric Spaces Involving Polynomial-Type Contractions}

\author[Gopinath Janardhanan ]{Gopinath Janardhanan$^{1}$}
\author[Gunaseelan Mani]{Gunaseelan Mani$^{1}$}
\author[Nancy Delaila John Kennedy]{Nancy Delaila John Kennedy$^{2}$}
\author[Ya\'e Ulrich Gaba]{Ya\'e Ulrich Gaba$^{3,4,\dagger}$}

	\address{$^{1}$ Department of Mathematics, \\Saveetha School of Engineering,\\ Saveetha Institute of Medical and Technical Sciences, SIMATS,\\ Chennai 602105, Tamilnadu, India.}

    \address{$^{2}$ Department of Mathematics, \\St. Joseph's College of Engineering \\OMR, Chennai-600119, Tamilnadu, India}
	
	\address{$^{3}$ AI Research and Innovation Nexus for Africa (AIRINA Labs) AI.Technipreneurs, Cotonou, B\'enin}

\address{$^{4}$ African Center for Advanced Studies (ACAS),		P.O. Box 4477, Yaoundé, Cameroon.}

    \email{\textcolor[rgb]{0.00,0.00,0.84}{gopijana2055@gmail.com
	}}
	
	\email{\textcolor[rgb]{0.00,0.00,0.84}{mathsguna@yahoo.com}}

    \email{\textcolor[rgb]{0.00,0.00,0.84}{nancydelaila20@gmail.com}}
	\email{\textcolor[rgb]{0.00,0.00,0.84}{yaeulrich.gaba@gmail.com
	}}

	
	
	\date{Received: xxxxxx; Accepted: zzzzzz.
		\newline \indent $^{\dagger}$Corresponding author}
	
	\begin{abstract}
    In this paper, we investigate the existence and uniqueness of fixed points for self-mappings defined on bipolar metric spaces using a new class of contractive conditions, namely polynomial-type contractions. Our main results establish sufficient conditions under which a mapping on a complete bipolar metric space admits a UFP. Several illustrative examples are provided to demonstrate the applicability of our theorems, and we further show how our results generalize and improve upon existing fixed point theorems in both standard and generalized metric settings.
\end{abstract}

	\maketitle

\section{Introduction}
The most commonly used techniques for establishing the existence and uniqueness of solutions to nonlinear problems—such as differential equations, integral equations, evolution equations—typically and fractional differential equations involve reducing the problem to an equation of the form $\breve{F}\mathit{e} = \mathit{e}$ and $\mathit{e} \in \breve{E}$ is the unknown solution.This is referred to as a FP of $\breve{F}$. \\

The Banach Fixed Point Theorem (FPT) stands as one of the most fundamental and widely recognized results in FP theory \cite{02}, which states that, if $(\breve{E}, \vartheta)$ is a complete $\mathcal{MS}$ and $\breve{F}\colon\breve{E}\rightarrow\breve{E}$ is a mapping satisfies
\begin{align*}
	\vartheta(\breve{F}\mathit{e}, \breve{F}\mathit{f})\leq\pi\vartheta(\mathit{e}, \mathit{f})
\end{align*}
for all $\mathit{e}, \mathit{f}\in\breve{E}$ and $\pi\in(0,1)$ is a constant, then
\begin{enumerate}[label=($M_\arabic*)$]
	\item $\breve{F}$ has a UFP;
	\item For all $\mathit{e}_0\in\breve{E}$, the sequence $\{\mathit{e}_\kappa\}$ given by $\mathit{e}_{\kappa+1}=\breve{F}\mathit{e}_\kappa$, approaches to FP.
\end{enumerate}
Agrawal et. al \cite{01} has established FPT in $\mathcal{MS}$. Banach \cite{02} has presented banach contraction principle. Berinde \cite{04} has established general constructive FPTs in \'{C}iri\'{c}-type almost contractions in $\mathcal{MS}$s. Berinde \cite{05} has proved FPT in banach
spaces by using a retraction-displacement condition. Berinde et. al \cite{06} has proved FPTs in almost contractions. Berinde et. al \cite{07} has proved FPs of enriched contractions in banach spaces. Boyd et. al \cite{08} has proved FPT in metrically complete $\mathcal{MS}$.

The concept of FPs has been a central theme in nonlinear analysis and has led to significant developments in various branches of mathematics. In 2000, Branciari \cite{09} introduced a FPT of Banach-Caccioppoli type in a class of generalized $\mathcal{MS}$s, which extends the classical Banach contraction principle. Earlier, Chatterjea \cite{10} presented a FP result based on a distinct type of contraction condition, now known as the Chatterjea contraction. In a further generalization, \v{C}iri\'{c} \cite{11} proposed a broader class of contractive mappings that unify and extend several known FPTs. The structure of generalized $\mathcal{MS}$s was further enriched by Czerwik \cite{12}, who explored contraction mappings in \( b \)-$\mathcal{MS}$s, thereby opening new avenues for FP theory. Additionally, Dhage \cite{13} introduced and analyzed mappings in generalized $\mathcal{MS}$s, contributing valuable FP results applicable in more abstract settings. These foundational works collectively form the basis for many contemporary studies in FP theory and its applications.

FPT continues to evolve with the development of novel contraction principles and generalizations applicable to diverse mathematical structures. Petrov \cite{14} introduced a geometric approach by studying mappings that contract the perimeters of triangles, offering a unique perspective within metric frameworks. Petru\c{s}el and Rus \cite{15} expanded FPT by incorporating both metric and order structures, enabling broader applicability in ordered $\mathcal{MS}$s. Popescu and P\u{a}curar \cite{16} recently proposed FP results for generalized Chatterjea-type mappings, contributing to ongoing efforts to unify and generalize contraction conditions. Classical results by Rakotch \cite{17} and Reich \cite{18} laid essential groundwork in the theory of contractive functions, which continues to inspire modern generalizations.

The idea of weakly Picard mappings was explored by Rus \cite{19}, leading to a series of studies on Picard operators \cite{20, 21, 22} that examined their properties, equivalences, and applications. These contributions underscore the importance of iterative behavior and convergence in FP theory.  Recent advances in bipolar metric spaces have extended classical fixed point results to accommodate asymmetric structures. Gaba et al.~\cite{g1} established Banach-type fixed point theorems in this context, while Aphane, and et al.~\cite{g2} introduced $(\alpha, BK)$-contractions, broadening the range of applicable contractive conditions. Zhang and Song \cite{23} further enriched the field by introducing the framework of generalized $(\psi-\phi)$-weak contractions, which accommodate a wider class of nonlinear mappings. Together, these works highlight the depth and versatility of FP theory in both classical and modern mathematical analysis. Gunaseelan mani et. al \cite{24} has presented FPT in $F$-contraction in bipolar $\mathcal{MS}$. Gunaseelan mani et. al \cite{25} has established FPT in controlled bipolar $\mathcal{MS}$. Gunaseelan mani et. al. \cite{26} has given common FPT in bipolar orthogonal $\mathcal{MS}$. Further more details see [\cite{{27},{28},{29},{30},{31}}]. Throughout in this paper $\mathcal{MS}$ means metric space and PC means polynomial contraction.   
\section{Preliminaries}
\begin{defn}
	Let $\breve{E}$ and $\breve{P}$ be a two non-void set and $\vartheta\colon\breve{E}\times\breve{P}\rightarrow[0, +\infty)$ be a mapping implies that
\begin{enumerate}[label=(\arabic*)]
	\item $\mathit{e}=\mathit{f}$ iff $\vartheta(\mathit{e}, \mathit{f})=0$;
	\item $\vartheta(\mathit{e}, \mathit{f})=\vartheta(\mathit{f}, \mathit{e})$ if $\mathit{e}, \mathit{f}\in\breve{E}\cap\breve{P}$;
	\item $\vartheta(\mathit{e}, \mathit{f})\leq\vartheta(\mathit{e}, \mathfrak{z})+\vartheta(\mathfrak{r}, \mathfrak{z})+\vartheta(\mathfrak{r}, \mathit{f})$, $\forall\mathit{e},\mathfrak{r}\in\breve{E}$ and $\mathfrak{z}, \mathit{f}\in\breve{P}$.
\end{enumerate}
Then, the mapping $\vartheta$ is called a bipolar $\mathcal{MS}$ on the pair $(\breve{E}, \breve{P})$ and the triple $(\breve{E}, \breve{P}, \vartheta)$ is called bipolar $\mathcal{MS}$.
\end{defn}
\begin{defn}
	Let $(\breve{E}, \vartheta)$ be a $\mathcal{MS}$. A mapping $\breve{F}\colon\breve{E}\rightarrow\breve{E}$ is called an almost contraction, if there exist $\pi\in(0,1)$ and $\rho>0$ implies that
\begin{align*}
	\vartheta(\breve{F}\mathit{e}, \breve{F}\mathit{f})\leq\pi\vartheta(\mathit{e}, \mathit{f})+\rho\vartheta(\mathit{f}, \breve{F}\mathit{e}), \forall\mathit{e}, \mathit{f}\in\breve{E}.
\end{align*}
\end{defn}
The following FPT for the above class of mappings was proven by Berinde \cite{03}.
\begin{thm}
	Let $(\breve{E}, \vartheta)$ be a complete $\mathcal{MS}$ and the mapping $\breve{F}\colon\breve{E}\rightarrow\breve{E}$ be an almost contraction. Then
\begin{enumerate}[label=(\roman*)]
	\item $\breve{F}$ admits at least one FP;
	\item For all $\{\mathit{g}_\kappa\}$ defined by $\mathit{g}_{\kappa+1}=\breve{F}\mathit{g}_\kappa$, converges to a FP of $\breve{F}$.
\end{enumerate}
\end{thm}

\section{The class of polynomial contractions(PC)}
 \begin{defn}
 	Let $(\breve{E}, \vartheta)$ be a bipolar $\mathcal{MS}$ and $\breve{F}\colon\breve{E}\cup\breve{P}\rightarrow\breve{E}\cup\breve{P}$ be a mapping. It is stated that $\breve{F}$ is a PC, if there exist $\pi\in(0, 1)$, $\sigma\geq1$ and a mappings $\mathit{q}_\upsilon\colon\breve{E}\times\breve{P}\rightarrow[0, \infty), \upsilon=0, \cdots\sigma$ implies that
\begin{align*}
	\sum_{\upsilon=0}^{\sigma}\mathit{q}_\upsilon(\breve{F}\mathit{e}, \breve{F}\mathit{f})\vartheta^\upsilon(\breve{F}\mathit{e}, \breve{F}\mathit{f})\leq\pi\sum_{\upsilon=0}^{\sigma}\mathit{q}_\upsilon(\mathit{e}, \mathit{f})\vartheta^\upsilon(\mathit{e}, \mathit{f})
\end{align*}
forall $\mathit{e}\in\breve{E}, \mathit{f}\in\breve{P}$.
 \end{defn} 
In this section, we focus on the study of FPs within bipolar $\mathcal{MS}$ involving the class of polynomial functions.
\begin{thm} \label{thm3.1}
	Let $(\breve{E}, \breve{P}, \vartheta)$ be a complete bipolar $\mathcal{MS}$ and the covariant mapping $\breve{F}\colon\breve{E}\cup\breve{P}\rightrightarrows\breve{E}\cup\breve{P}$ be a PC such that
\begin{enumerate}[label=(\roman*)]
	\item $\breve{F}$ is continuous;
	\item We can find that $\varrho\in\{1,\cdots,\sigma\}$ and $\breve{Q}_\varrho>0$ implies that \label{thm3c2}
\begin{align*}
	\mathit{q}_\varrho(\mathit{e}, \mathit{f})\geq\breve{Q}_\varrho, \mathit{e}\in\breve{E}, \mathit{f}\in\breve{P}.
\end{align*}
\end{enumerate}
Then, $\breve{F}$ admits a UFP. Moreover for every, $\mathit{g}_0\in\breve{E}$, the picard sequence $\{\mathit{g}_\kappa\}\subset\breve{E}$ by $\mathit{g}_{\kappa+1}=\breve{F}\mathit{g}_\kappa$ and $\mathit{h}_0\in\breve{P}$, the picard sequence $\{\mathit{h}_\kappa\}\subset\breve{P}$ by $\mathit{h}_{\kappa+1}=\breve{F}\mathit{h}_\kappa, \forall\kappa\geq0$.
\begin{proof}
Initially, we show that the set of FPs of $\breve{F}$ is non-empty. Let $\mathit{g}_0\in\breve{E}$ and $\mathit{h}_0\in\breve{P}$ be FPs and $\{\mathit{g}_\kappa\}\subset\breve{E}$ and $\{\mathit{h}_\kappa\}\subset\breve{P}$ is defined by
\begin{align*}
	\mathit{g}_{\kappa+1}=\breve{F}\mathit{g}_\kappa \quad\text{and}\quad\mathit{h}_{\kappa+1}=\breve{F}\mathit{h}_\kappa, \quad\kappa\geq0.
\end{align*}
By definition of PC, with $(\mathit{e}, \mathit{f})=(\mathit{g}_0, \mathit{h}_0)$, we get
\begin{align*}
	\sum_{\upsilon=0}^{\sigma}\mathit{q}_\upsilon(\breve{F}\mathit{g}_0, \breve{F}\mathit{h}_0)\vartheta^\upsilon(\breve{F}\mathit{g}_0, \breve{F}\mathit{h}_0)\leq\pi\sum_{\upsilon=0}^{\sigma}\mathit{q}_\upsilon(\mathit{g}_0, \mathit{h}_0)\vartheta^\upsilon(\mathit{g}_0, \mathit{h}_0),
\end{align*}
i.e.,
\begin{align*}
	\sum_{\upsilon=0}^{\sigma}\mathit{q}_\upsilon(\mathit{g}_1, \mathit{h}_1)\vartheta^\upsilon(\mathit{g}_1, \mathit{h}_1)\leq\pi\sum_{\upsilon=0}^{\sigma}\mathit{q}_\upsilon(\mathit{g}_0, \mathit{h}_0)\vartheta^\upsilon(\mathit{g}_0, \mathit{h}_0).
\end{align*}
Generally we can write,
\begin{align*}
	\vartheta^\varrho(\mathit{g}_\kappa, \mathit{h}_{\kappa})\leq\pi^\kappa\sum_{\upsilon=0}^{\sigma}\mathit{q}_\upsilon(\mathit{g}_0, \mathit{h}_0)\vartheta^\upsilon(\mathit{g}_0, \mathit{h}_0), \quad\kappa\geq0.
\end{align*}
By definition of PC, with $(\mathit{e}, \mathit{f})=(\mathit{g}_0, \mathit{h}_1)$, we get
\begin{align*}
	\sum_{\upsilon=0}^{\sigma}\mathit{q}_\upsilon(\breve{F}\mathit{g}_0, \breve{F}\mathit{h}_1)\vartheta^\upsilon(\breve{F}\mathit{g}_0, \breve{F}\mathit{h}_1)\leq\pi\sum_{\upsilon=0}^{\sigma}\mathit{q}_\upsilon(\mathit{g}_0, \mathit{h}_1)\vartheta^\upsilon(\mathit{g}_0, \mathit{h}_1),
\end{align*}
i.e.,
\begin{align}
	\sum_{\upsilon=0}^{\sigma}\mathit{q}_\upsilon(\mathit{g}_1, \mathit{h}_2)\vartheta^\upsilon(\mathit{g}_1, \mathit{h}_2)\leq\pi\sum_{\upsilon=0}^{\sigma}\mathit{q}_\upsilon(\mathit{g}_0, \mathit{h}_1)\vartheta^\upsilon(\mathit{g}_0, \mathit{h}_1). \label{eq1}
\end{align}
Again by definition of PC with $(\mathit{e}, \mathit{f})=(\mathit{g}_1, \mathit{h}_2)$
\begin{align*}
	\sum_{\upsilon=0}^{\sigma}\mathit{q}_\upsilon(\breve{F}\mathit{g}_1, \breve{F}\mathit{h}_2)\vartheta^\upsilon(\breve{F}\mathit{g}_1, \breve{F}\mathit{h}_2)\leq\pi\sum_{\upsilon=0}^{\sigma}\mathit{q}_\upsilon(\mathit{g}_1, \mathit{h}_2)\vartheta^\upsilon(\mathit{g}_1, \mathit{h}_2),
\end{align*}
i.e., 
\begin{align*}
	\sum_{\upsilon=0}^{\sigma}\mathit{q}_\upsilon(\mathit{g}_2, \mathit{h}_3)\vartheta^\upsilon(\mathit{g}_2, \mathit{h}_3)\leq\pi\sum_{\upsilon=0}^{\sigma}\mathit{q}_\upsilon(\mathit{g}_1, \mathit{h}_2)\vartheta^\upsilon(\mathit{g}_1, \mathit{h}_2),
\end{align*}
which implies
\begin{align*}
	\sum_{\upsilon=0}^{\sigma}\mathit{q}_\upsilon(\mathit{g}_2, \mathit{h}_3)\vartheta^\upsilon(\mathit{g}_2, \mathit{h}_3)\leq\pi^2\sum_{\upsilon=0}^{\sigma}\mathit{q}_\upsilon(\mathit{g}_0, \mathit{h}_1)\vartheta^\upsilon(\mathit{g}_0, \mathit{h}_1).
\end{align*}
Continuing in the same way, we get
\begin{align}
	\sum_{\upsilon=0}^{\sigma}\mathit{q}_\upsilon(\mathit{g}_\kappa, \mathit{h}_{\kappa+1})\vartheta^\upsilon(\mathit{g}_\kappa, \mathit{h}_{\kappa+1})\leq\pi^\kappa\sum_{\upsilon=0}^{\sigma}\mathit{q}_\upsilon(\mathit{g}_0, \mathit{h}_1)\vartheta^\upsilon(\mathit{g}_0, \mathit{h}_1). \label{eq2}
\end{align}
Since
\begin{align*}
	\mathit{q}_\varrho(\mathit{g}_\kappa, \mathit{h}_{\kappa+1})\vartheta^\varrho(\mathit{g}_\kappa, \mathit{h}_{\kappa+1})\leq\sum_{\upsilon=0}^{\sigma}\mathit{q}_\upsilon(\mathit{g}_\kappa, \mathit{h}_{\kappa+1})\vartheta^\upsilon(\mathit{g}_\kappa, \mathit{h}_{\kappa+1}). 
\end{align*}
We obtain by \ref{thm3c2} that
\begin{align*}
	\breve{Q}_\varrho\vartheta^\varrho(\mathit{g}_\kappa, \mathit{h}_{\kappa+1})\leq\sum_{\upsilon=0}^{\sigma}\mathit{q}_\upsilon(\mathit{g}_\kappa, \mathit{h}_{\kappa+1})\vartheta^\upsilon(\mathit{g}_\kappa, \mathit{h}_{\kappa+1}),
\end{align*}
which implies by \eqref{eq2}, that
\begin{align}
	\vartheta^\varrho(\mathit{g}_\kappa, \mathit{h}_{\kappa+1})\leq\pi^\kappa\sum_{\upsilon=0}^{\sigma}\mathit{q}_\upsilon(\mathit{g}_0, \mathit{h}_1)\vartheta^\upsilon(\mathit{g}_0, \mathit{h}_1), \quad\kappa\geq0. \label{eq3}
\end{align}
Similarly we have,
\begin{align*}
	\vartheta^\varrho(\mathit{g}_\kappa, \mathit{h}_{\kappa})\leq\pi^\kappa\sum_{\upsilon=0}^{\sigma}\mathit{q}_\upsilon(\mathit{g}_0, \mathit{h}_0)\vartheta^\upsilon(\mathit{g}_0, \mathit{h}_0), \quad\kappa\geq0,
\end{align*}
where,
\begin{align}
	\mathcal{M}=\breve{Q}^{-1}_\varrho\sum_{\upsilon=0}^{\sigma}\mathit{q}_\upsilon(\mathit{g}_0, \mathit{h}_1)\vartheta^\upsilon(\mathit{g}_0, \mathit{h}_1)+\breve{Q}^{-1}_\varrho\sum_{\upsilon=0}^{\sigma}\mathit{q}_\upsilon(\mathit{g}_0, \mathit{h}_0)\vartheta^\upsilon(\mathit{g}_0, \mathit{h}_0). \label{eq4}
\end{align}
Then, by \eqref{eq3} and triangle inequality, we get
\begin{align*}
	\vartheta^\varrho(\mathit{g}_{\kappa+\varpi}, \mathit{h}_\kappa)&\leq\vartheta^\varrho(\mathit{g}_{\kappa+\varpi}, \mathit{h}_{\kappa+1})+\vartheta^\varrho(\mathit{g}_\kappa, \mathit{h}_{\kappa+1})+\vartheta^\varrho(\mathit{g}_\kappa, \mathit{h}_\kappa)\\
	&\leq\vartheta^\varrho(\mathit{g}_{\kappa+\varpi}, \mathit{h}_{\kappa+1})+\pi^\kappa\mathcal{M}\\
	&\leq\vartheta^\varrho(\mathit{g}_{\kappa+\varpi}, \mathit{h}_{\kappa+2})+\vartheta^\varrho(\mathit{g}_{\kappa+1}, \mathit{h}_{\kappa+2})+\vartheta^\varrho(\mathit{g}_{\kappa+1}, \mathit{h}_{\kappa+1})+\pi^\kappa\mathcal{M}\\
	&\leq\vartheta^\varrho(\mathit{g}_{\kappa+\varpi}, \mathit{h}_{\kappa+2})+(\pi^{\kappa+1}+\pi^\kappa)\mathcal{M}\\
	& \vdots \\
	&\leq\vartheta^\varrho(\mathit{g}_{\kappa+\varpi}, \mathit{h}_{\kappa+\varpi})+(\pi^{\kappa+\varpi-1}+\cdots+\pi^{\kappa+1}+\pi^\kappa)\mathcal{M}\\
	&\leq(\pi^{\kappa+\varpi}+\cdots+\pi^{\kappa+1}+\pi^\kappa)\mathcal{M}\\
	&\leq\mathcal{M}\frac{\pi^\kappa}{1-\pi},
\end{align*}
which implies,
\begin{align*}
	\vartheta(\mathit{g}_{\kappa+\varpi}, \mathit{h}_{\kappa})\leq\bigg(\frac{\mathcal{M}}{1-\pi}\bigg)^\frac{1}{\varrho}\bigg(\pi^\frac{1}{\varrho}\bigg)^\kappa\rightarrow0 \quad\text{as}\quad\kappa,\varpi\rightarrow\infty.
\end{align*}
Similarly,
\begin{align*}
	\vartheta(\mathit{g}_{\kappa}, \mathit{h}_{\kappa+\varpi})\leq\bigg(\frac{\mathcal{M}}{1-\pi}\bigg)^\frac{1}{\varrho}\bigg(\pi^\frac{1}{\varrho}\bigg)^\kappa\rightarrow0 \quad\text{as}\quad\kappa,\varpi\rightarrow\infty.
\end{align*}
This shows that $\{\mathit{g}_\kappa\}$ and $\{\mathit{h}_\kappa\}$ is a Cauchy bisequence. Since $(\breve{E}, \breve{P}, \vartheta)$ is complete, there exist $\mathit{g}\in\breve{E}$ such that
\begin{align*}
	\lim_{\kappa\rightarrow\infty}\vartheta(\mathit{g}_\kappa, \mathit{g})=0,
\end{align*}
due to the continuity of $\breve{F}$ that
\begin{align*}
	\lim_{\kappa\rightarrow\infty}\vartheta(\mathit{g}_{\kappa+1}, \breve{F}\mathit{g})=\lim_{\kappa\rightarrow\infty}\vartheta(\breve{F}\mathit{g}_\kappa, \breve{F}\mathit{g})=0.
\end{align*}
By the uniqueness of the limit, we conclude that $\breve{F}\mathit{g}=\mathit{g}$ i.e., $\mathit{g}$ is a FP of $\breve{F}$.\\
Now, we show that $\mathit{g}$ is the UFP of $\breve{F}$. Indeed, if $\mathit{g}^*\in\breve{E}$ is another FP of $\breve{F}$ i.e., $\breve{F}\mathit{g}^*=\mathit{g}^*$ and $\vartheta(\mathit{g}, \mathit{g}^*)>0$ then by definition of PC with $(\mathit{e}, \mathit{f})=(\mathit{g}, \mathit{g}^*)$
\begin{align*}
	\sum_{\upsilon=0}^{\sigma}\mathit{q}_\upsilon(\breve{F}\mathit{g}, \breve{F}\mathit{g}^*)\vartheta^\upsilon(\breve{F}\mathit{g}, \breve{F}\mathit{g}^*)\leq\pi\sum_{\upsilon=0}^{\sigma}\mathit{q}_\upsilon(\mathit{g}, \mathit{g}^*)\vartheta^\upsilon(\mathit{g}, \mathit{g}^*),
\end{align*}
i.e.,
\begin{align}
	\sum_{\upsilon=0}^{\sigma}\mathit{q}_\upsilon(\mathit{g}, \mathit{g}^*)\vartheta^\upsilon(\mathit{g}, \mathit{g}^*)\leq\pi\sum_{\upsilon=0}^{\sigma}\mathit{q}_\upsilon(\mathit{g}, \mathit{g}^*)\vartheta^\upsilon(\mathit{g}, \mathit{g}^*). \label{eq5}
\end{align}
On the other hand, from \ref{thm3c2}, we get
\begin{align*}
	\sum_{\upsilon=0}^{\sigma}\mathit{q}_\upsilon(\mathit{g}, \mathit{g}^*)\vartheta^\upsilon(\mathit{g}, \mathit{g}^*)&\geq\mathit{q}_\varrho(\mathit{g}, \mathit{g}^*)\vartheta^\upsilon(\mathit{g}, \mathit{g}^*)\\
	&\geq\breve{Q}_\varrho\vartheta^\upsilon(\mathit{g}, \mathit{g}^*).
\end{align*}
Since $\breve{Q}_\varrho>0$ and $\vartheta(\mathit{g}, \mathit{g}^*)>0$, we deduce that
\begin{align*}
	\sum_{\upsilon=0}^{\sigma}\mathit{q}_\upsilon(\mathit{g}, \mathit{g}^*)\vartheta^\upsilon(\mathit{g}, \mathit{g}^*)>0.
\end{align*}
Then, dividing by \eqref{eq5}, we get $\sum_{\upsilon=0}^{\sigma}\mathit{q}_\upsilon(\mathit{g}, \mathit{g}^*)\vartheta^\upsilon(\mathit{g}, \mathit{g}^*)$, we reach a contraction with $\pi\in(0, 1)$. Consequently, $\mathit{g}$ is a UFP of $\breve{F}$. Hence its completes.
\end{proof}
\end{thm}
\begin{thm} 
	Let $(\breve{E}, \breve{P}, \vartheta)$ be a complete bipolar $\mathcal{MS}$ and the contravariant mapping $\breve{F}\colon\breve{E}\cup\breve{P}\leftrightarrows\breve{E}\cup\breve{P}$ be a PC such that
	\begin{enumerate}[label=(\roman*)]
		\item $\breve{F}$ is continuous;
		\item We can find that $\varrho\in\{1,\cdots,\sigma\}$ and $\breve{Q}_\varrho>0$ implies that \label{thm3.2c2}
		\begin{align*}
			\mathit{q}_\varrho(\mathit{e}, \mathit{f})\geq\breve{Q}_\varrho, \mathit{e}\in\breve{E}, \mathit{f}\in\breve{P}.
		\end{align*}
	\end{enumerate}
	Then, $\breve{F}$ admits a UFP. Moreover for every, $\mathit{g}_0\in\breve{E}$, the picard sequence $\{\mathit{g}_\kappa\}\subset\breve{E}$ by $\mathit{g}_{\kappa+1}=\breve{F}\mathit{g}_\kappa$ and $\mathit{h}_0\in\breve{P}$, the picard sequence $\{\mathit{h}_\kappa\}\subset\breve{P}$ by $\mathit{h}_{\kappa+1}=\breve{F}\mathit{h}_\kappa, \forall\kappa\geq0$.
	\begin{proof}
		Initially, we show that the set of FPs of $\breve{F}$ is non-empty. Let $\mathit{g}_0\in\breve{E}$ and $\mathit{h}_0\in\breve{P}$ be FPs and $\{\mathit{g}_\kappa\}\subset\breve{E}$ and $\{\mathit{h}_\kappa\}\subset\breve{P}$ is defined by
		\begin{align*}
			\mathit{h}_{\kappa}=\breve{F}\mathit{g}_\kappa \quad\text{and}\quad\mathit{g}_{\kappa+1}=\breve{F}\mathit{h}_\kappa, \quad\kappa\geq0.
		\end{align*}
		By definition of PC, with $(\mathit{e}, \mathit{f})=(\mathit{g}_0, \mathit{h}_0)$, we get
		\begin{align*}
			\sum_{\upsilon=0}^{\sigma}\mathit{q}_\upsilon(\breve{F}\mathit{g}_0, \breve{F}\mathit{h}_0)\vartheta^\upsilon(\breve{F}\mathit{g}_0, \breve{F}\mathit{h}_0)\leq\pi\sum_{\upsilon=0}^{\sigma}\mathit{q}_\upsilon(\mathit{g}_0, \mathit{h}_0)\vartheta^\upsilon(\mathit{g}_0, \mathit{h}_0),
		\end{align*}
		i.e.,
		\begin{align}
			\sum_{\upsilon=0}^{\sigma}\mathit{q}_\upsilon(\mathit{h}_0, \mathit{g}_1)\vartheta^\upsilon(\mathit{h}_0, \mathit{g}_1)\leq\pi\sum_{\upsilon=0}^{\sigma}\mathit{q}_\upsilon(\mathit{g}_0, \mathit{h}_0)\vartheta^\upsilon(\mathit{g}_0, \mathit{h}_0). \label{eq6}
		\end{align}
Generally we can write,
\begin{align*}
	\vartheta^\varrho(\mathit{g}_{\kappa}, \mathit{h}_{\kappa})\leq\pi^{2\kappa}\sum_{\upsilon=0}^{\sigma}\mathit{q}_\upsilon(\mathit{g}_0, \mathit{h}_1)\vartheta^\upsilon(\mathit{g}_0, \mathit{h}_1), \quad\kappa\geq0. 
\end{align*}
		Again by definition of PC with $(\mathit{e}, \mathit{f})=(\mathit{g}_0, \mathit{h}_1)$
		\begin{align*}
			\sum_{\upsilon=0}^{\sigma}\mathit{q}_\upsilon(\breve{F}\mathit{g}_0, \breve{F}\mathit{h}_1)\vartheta^\upsilon(\breve{F}\mathit{g}_0, \breve{F}\mathit{h}_1)\leq\pi\sum_{\upsilon=0}^{\sigma}\mathit{q}_\upsilon(\mathit{g}_0, \mathit{h}_1)\vartheta^\upsilon(\mathit{g}_0, \mathit{h}_1),
		\end{align*}
		i.e., 
		\begin{align*}
			\sum_{\upsilon=0}^{\sigma}\mathit{q}_\upsilon(\mathit{h}_0, \mathit{g}_2)\vartheta^\upsilon(\mathit{h}_0, \mathit{g}_2)\leq\pi\sum_{\upsilon=0}^{\sigma}\mathit{q}_\upsilon(\mathit{g}_0, \mathit{h}_1)\vartheta^\upsilon(\mathit{g}_0, \mathit{h}_1),
		\end{align*}
		which implies
		\begin{align*}
			\sum_{\upsilon=0}^{\sigma}\mathit{q}_\upsilon(\mathit{g}_2, \mathit{h}_0)\vartheta^\upsilon(\mathit{g}_2, \mathit{h}_0)\leq\pi^2\sum_{\upsilon=0}^{\sigma}\mathit{q}_\upsilon(\mathit{g}_0, \mathit{h}_1)\vartheta^\upsilon(\mathit{g}_0, \mathit{h}_1).
		\end{align*}
		Continuing in the same way, we get
		\begin{align}
			\sum_{\upsilon=0}^{\sigma}\mathit{q}_\upsilon(\mathit{g}_{\kappa+1}, \mathit{h}_{\kappa})\vartheta^\upsilon(\mathit{g}_{\kappa+1}, \mathit{h}_{\kappa})\leq\pi^{2\kappa+1}\sum_{\upsilon=0}^{\sigma}\mathit{q}_\upsilon(\mathit{g}_0, \mathit{h}_1)\vartheta^\upsilon(\mathit{g}_0, \mathit{h}_1). \label{eq7}
		\end{align}
		Since
		\begin{align*}
			\mathit{q}_\varrho(\mathit{g}_{\kappa+1}, \mathit{h}_{\kappa})\vartheta^\varrho(\mathit{g}_{\kappa+1}, \mathit{h}_{\kappa})\leq\sum_{\upsilon=0}^{\sigma}\mathit{q}_\upsilon(\mathit{g}_{\kappa+1}, \mathit{h}_{\kappa})\vartheta^\upsilon(\mathit{g}_{\kappa+1}, \mathit{h}_{\kappa}). 
		\end{align*}
		We obtain by \ref{thm3.2c2} that
		\begin{align*}
			\breve{Q}_\varrho\vartheta^\varrho(\mathit{g}_{\kappa+1}, \mathit{h}_{\kappa})\leq\sum_{\upsilon=0}^{\sigma}\mathit{q}_\upsilon(\mathit{g}_{\kappa+1}, \mathit{h}_{\kappa})\vartheta^\upsilon(\mathit{g}_{\kappa+1}, \mathit{h}_{\kappa}),
		\end{align*}
		which implies by \eqref{eq7}, that
		\begin{align}
			\vartheta^\varrho(\mathit{g}_{\kappa+1}, \mathit{h}_{\kappa})\leq\pi^{2\kappa+1}\sum_{\upsilon=0}^{\sigma}\mathit{q}_\upsilon(\mathit{g}_0, \mathit{h}_1)\vartheta^\upsilon(\mathit{g}_0, \mathit{h}_1), \quad\kappa\geq0,  \label{eq8}
		\end{align}
similarly we have,
\begin{align*}
	\vartheta^\varrho(\mathit{g}_{\kappa}, \mathit{h}_{\kappa})\leq\pi^{2\kappa}\sum_{\upsilon=0}^{\sigma}\mathit{q}_\upsilon(\mathit{g}_0, \mathit{h}_1)\vartheta^\upsilon(\mathit{g}_0, \mathit{h}_1), \quad\kappa\geq0,  
\end{align*}
		where,
		\begin{align}
			\mathcal{M}=\breve{Q}^{-1}_\varrho\sum_{\upsilon=0}^{\sigma}\mathit{q}_\upsilon(\mathit{g}_0, \mathit{h}_1)\vartheta^\upsilon(\mathit{g}_0, \mathit{h}_1)+\breve{Q}^{-1}_\varrho\sum_{\upsilon=0}^{\sigma}\mathit{q}_\upsilon(\mathit{g}_0, \mathit{h}_1)\vartheta^\upsilon(\mathit{g}_0, \mathit{h}_1). \label{eq9}
		\end{align}
		Then, by \eqref{eq8} and triangle inequality, we get
\begin{align*}
	\vartheta^\varrho(\mathit{g}_{\kappa+\varpi}, \mathit{h}_\kappa)&\leq\vartheta^\varrho(\mathit{g}_{\kappa+\varpi}, \mathit{h}_{\kappa+1})+\vartheta^\varrho(\mathit{g}_{\kappa+1}, \mathit{h}_{\kappa+1})+\vartheta^\varrho(\mathit{g}_{\kappa+1}, \mathit{h}_\kappa)\\
	&\leq\vartheta^\varrho(\mathit{g}_{\kappa+\varpi}, \mathit{h}_{\kappa+1})+(\pi^{2\kappa+2}+\pi^{2\kappa+1})\mathcal{M}\\
	&\leq\vartheta^\varrho(\mathit{g}_{\kappa+\varpi}, \mathit{h}_{\kappa+2})+\vartheta^\varrho(\mathit{g}_{\kappa+2}, \mathit{h}_{\kappa+2})+\vartheta^\varrho(\mathit{g}_{\kappa+2}, \mathit{h}_{\kappa+1})\\
	&+(\pi^{2\kappa+2}+\pi^{2\kappa+1})\mathcal{M}\\
	&\vdots\\
	&\leq\vartheta^\varrho(\mathit{g}_{\kappa+\varpi}, \mathit{h}_{\kappa+\varpi-1})+(\pi^{2\kappa+2\varpi-2}+\cdots+\pi^{2\kappa+1}+\pi^\kappa)\mathcal{M}\\
	&\leq(\pi^{2\kappa+2\varpi-1}+\pi^{2\kappa+2\varpi-2}+\pi^{2\kappa+2\varpi-3}+\cdots+\pi^{2\kappa+1})\mathcal{M}\\
	&\leq\mathcal{M}\frac{\pi^{2\kappa}}{1-\pi},
\end{align*}
		which implies,
		\begin{align*}
			\vartheta(\mathit{g}_{\kappa+\varpi}, \mathit{h}_{\kappa})\leq\bigg(\frac{\mathcal{M}}{1-\pi}\bigg)^\frac{1}{\varrho}\bigg(\pi^\frac{1}{\varrho}\bigg)^{2\kappa}\rightarrow0 \quad\text{as}\quad\kappa,\varpi\rightarrow\infty,
		\end{align*}
similarly,
\begin{align*}
	\vartheta(\mathit{g}_{\kappa}, \mathit{h}_{\kappa+\varpi})\leq\bigg(\frac{\mathcal{M}}{1-\pi}\bigg)^\frac{1}{\varrho}\bigg(\pi^\frac{1}{\varrho}\bigg)^{2\kappa}\rightarrow0 \quad\text{as}\quad\kappa,\varpi\rightarrow\infty.
\end{align*}
		This shows that $\{\mathit{g}_\kappa\}$ and $\{\mathit{h}_\kappa\}$ is a Cauchy bisequence. Since $(\breve{E}, \breve{P}, \vartheta)$ is complete, there exist $\mathit{g}\in\breve{E}$ such that
		\begin{align*}
			\lim_{\kappa\rightarrow\infty}\vartheta(\mathit{g}_\kappa, \mathit{g})=0,
		\end{align*}
		due to the continuity of $\breve{F}$ that
		\begin{align*}
			\lim_{\kappa\rightarrow\infty}\vartheta(\mathit{g}_{\kappa+1}, \breve{F}\mathit{g})=\lim_{\kappa\rightarrow\infty}\vartheta(\breve{F}\mathit{g}_\kappa, \breve{F}\mathit{g})=0.
		\end{align*}
		By the uniqueness of the limit, we conclude that $\breve{F}\mathit{g}=\mathit{g}$ i.e., $\mathit{g}$ is a FP of $\breve{F}$.\\
		Now, we show that $\mathit{g}$ is the UFP of $\breve{F}$. Indeed, if $\mathit{g}^*\in\breve{E}$ is another FP of $\breve{F}$ i.e., $\breve{F}\mathit{g}^*=\mathit{g}^*$ and $\vartheta(\mathit{g}, \mathit{g}^*)>0$ then by definition of PC with $(\mathit{e}, \mathit{f})=(\mathit{g}, \mathit{g}^*)$
		\begin{align*}
			\sum_{\upsilon=0}^{\sigma}\mathit{q}_\upsilon(\breve{F}\mathit{g}, \breve{F}\mathit{g}^*)\vartheta^\upsilon(\breve{F}\mathit{g}, \breve{F}\mathit{g}^*)\leq\pi\sum_{\upsilon=0}^{\sigma}\mathit{q}_\upsilon(\mathit{g}, \mathit{g}^*)\vartheta^\upsilon(\mathit{g}, \mathit{g}^*),
		\end{align*}
		i.e.,
		\begin{align}
			\sum_{\upsilon=0}^{\sigma}\mathit{q}_\upsilon(\mathit{g}, \mathit{g}^*)\vartheta^\upsilon(\mathit{g}, \mathit{g}^*)\leq\pi\sum_{\upsilon=0}^{\sigma}\mathit{q}_\upsilon(\mathit{g}, \mathit{g}^*)\vartheta^\upsilon(\mathit{g}, \mathit{g}^*). \label{eq10}
		\end{align}
		On the other hand, from \ref{thm3.2c2}, we get
		\begin{align*}
			\sum_{\upsilon=0}^{\sigma}\mathit{q}_\upsilon(\mathit{g}, \mathit{g}^*)\vartheta^\upsilon(\mathit{g}, \mathit{g}^*)&\geq\mathit{q}_\varrho(\mathit{g}, \mathit{g}^*)\vartheta^\upsilon(\mathit{g}, \mathit{g}^*)\\
			&\geq\breve{Q}_\varrho\vartheta^\upsilon(\mathit{g}, \mathit{g}^*).
		\end{align*}
		Since $\breve{Q}_\varrho>0$ and $\vartheta(\mathit{g}, \mathit{g}^*)>0$, we deduce that
		\begin{align*}
			\sum_{\upsilon=0}^{\sigma}\mathit{q}_\upsilon(\mathit{g}, \mathit{g}^*)\vartheta^\upsilon(\mathit{g}, \mathit{g}^*)>0.
		\end{align*}
		Then, dividing by \eqref{eq10}, we get $\sum_{\upsilon=0}^{\sigma}\mathit{q}_\upsilon(\mathit{g}, \mathit{g}^*)\vartheta^\upsilon(\mathit{g}, \mathit{g}^*)$, we reach a contraction with $\pi\in(0, 1)$. Consequently, $\mathit{g}$ is a UFP of $\breve{F}$. Hence its completes.
	\end{proof}
\end{thm}
\begin{prop} \label{pro1}
	Let $(\breve{E}, \breve{P}, \vartheta)$ be a bipolar $\mathcal{MS}$ and the mapping $\breve{F}\colon\breve{E}\cup\breve{P}\rightarrow\breve{E}\cup\breve{P}$ be a PC such that
\begin{enumerate}[label=(\roman*)]
	\item $\mathit{q}_0\equiv0, \quad\text{i.e.,}\quad\mathit{q}_0(\mathit{e}, \mathit{f})=0, \forall\mathit{e}\in\breve{E}, \mathit{f}\in\breve{P}$; \label{proc1}
	\item Forall, $\upsilon\in\{1,\cdots,\sigma\}$, there exist $\breve{W}_\upsilon>0$ such that \label{proc2}
\begin{align*}
	\mathit{q}_\upsilon(\mathit{e}, \mathit{f})\leq\breve{W}_\upsilon, \mathit{e}\in\breve{E}, \mathit{f}\in\breve{P}
\end{align*}
\item There exist $\varrho\in\{1,\cdots,\sigma\}$ and $\breve{Q}_\varrho>0$ such that \label{proc3}
\begin{align*}
	\mathit{q}_\varrho(\mathit{e}, \mathit{f})\geq\breve{Q}_\varrho, \mathit{e}\in\breve{E}, \mathit{f}\in\breve{P}.
\end{align*}
\end{enumerate}
Then $\breve{F}$ is continuous.
\begin{proof}
	Let $\{\mathit{g}_\kappa\}\subset\breve{E}$ and $\{\mathit{h}_\kappa\}\subset\breve{P}$ be a bisequence such that
\begin{align}
	\lim_{\kappa\rightarrow\infty}\vartheta(\mathit{g}_\kappa, \mathit{h}_\kappa)=0, \label{eq11}
\end{align}
for some $\mathit{g}\in\breve{E}$ and $\mathit{h}\in\breve{P}$. Using \ref{proc1} and by definition of PC with $(\mathit{e}, \mathit{f})=(\mathit{g}_\kappa, \mathit{h}_\kappa)$, we get
\begin{align*}
	\sum_{\upsilon=0}^{\sigma}\mathit{q}_\upsilon(\breve{F}\mathit{g}_\kappa, \breve{F}\mathit{h}_\kappa)\vartheta^\upsilon(\breve{F}\mathit{g}_\kappa, \breve{F}\mathit{h}_\kappa)\leq\pi\sum_{\upsilon=0}^{\sigma}\mathit{q}_\upsilon(\mathit{g}_\kappa, \mathit{h}_\kappa)\vartheta^\upsilon(\mathit{g}_\kappa, \mathit{h}_\kappa), \kappa\geq0,
\end{align*}
 by \ref{proc2} and \ref{proc3} implies that
\begin{align}
	\breve{Q}_\varrho\vartheta^\varrho(\breve{F}\mathit{g}_\kappa,\breve{F}\mathit{h}_\kappa)\leq\pi\sum_{\upsilon=1}^{\sigma}\breve{W}_\upsilon\vartheta^\upsilon(\mathit{g}_\kappa, \mathit{h}_\kappa), \kappa\geq0. \label{eq12}
\end{align}
Then, by \eqref{eq11} and taking the limit as $\kappa\rightarrow\infty$ in \eqref{eq12}, we obtain
\begin{align*}
	\lim_{\kappa\rightarrow\infty}\vartheta^\varrho(\breve{F}\mathit{g}_\kappa, \breve{F}\mathit{h}_\kappa)=0,
\end{align*}
which is equivalent to,
\begin{align*}
	\lim_{\kappa\rightarrow\infty}\vartheta(\breve{F}\mathit{g}_\kappa, \breve{F}\mathit{h}_\kappa)=0.
\end{align*}
Thus $\breve{F}$ is a continuous mapping.
\end{proof}
\end{prop}
The following result is derived from Theorem \ref{thm3.1} and Proposition \ref{pro1}.
\begin{cor}
	Let $(\breve{E}, \breve{P}, \vartheta)$ be a complete bipolar $\mathcal{MS}$ and $\breve{F}\colon\breve{E}\cup\breve{P}\rightarrow\breve{E}\cup\breve{P}$ be a polynomial contaction such that
\begin{enumerate}[label=(\roman*)]
	\item $\mathit{q}_0\equiv0$;
	\item For all $\upsilon\in\{1,\cdots,\sigma\}$, there exists $\breve{W}_\upsilon>0$ such that
\begin{align*}
	\mathit{q}_\upsilon(\mathit{e}, \mathit{f})\leq\breve{W}_\upsilon, \mathit{e}\in\breve{E}, \mathit{f}\in\breve{P};
\end{align*}
\item There exist $\varrho\in\{1,\cdots,\sigma\}$ and $\breve{Q}_\varrho>0$ such that
\begin{align*}
	\mathit{q}_\varrho(\mathit{e}, \mathit{f})\leq\breve{Q}_\varrho, \mathit{e}\in\breve{E}, \mathit{f}\in\breve{P};
\end{align*}
\end{enumerate}
Then $\breve{F}$ admits a UFP.
\end{cor}
The following result follows directly from the above corollary.
\begin{cor}
	Let $(\breve{E}, \breve{P}, \vartheta)$ be a complete bipolar $\mathcal{MS}$ and $\breve{F}\colon\breve{E}\cup\breve{P}\rightarrow\breve{E}\cup\breve{P}$ be a mapping. Consider we can find that $\pi\in(0, 1)$,  $\sigma\geq1$ and a sequence $\{\mathit{q}_\upsilon\}_{\upsilon=1}^{\sigma}\subset(0, \infty)$ implies that
\begin{align*}
	\sum_{\upsilon=1}^{\sigma}\mathit{q}_\upsilon\vartheta^\upsilon(\breve{F}\mathit{e}, \breve{F}\mathit{f})\leq\pi\sum_{\upsilon=1}^{\sigma}\mathit{q}_\upsilon\vartheta^\upsilon(\mathit{e}, \mathit{f}), \forall\mathit{e}\in\breve{E}, \mathit{f}\in\breve{P}
\end{align*}
Then $\breve{F}$ admits a UFP.
\end{cor}
\begin{example}
	Let $\breve{E}=\{\mathit{e}_1, \mathit{e}_2, \mathit{e}_3\}$ and $\mathit{f}=\{\mathit{e}_1, \mathit{e}_5\}$ and the mapping be $\breve{F}\colon\breve{E}\cup\breve{P}\rightarrow\breve{E}\cup\breve{P}$ is defined by
\begin{align*}
	\breve{F}\mathit{e}_1=\mathit{e}_1, \breve{F}\mathit{e}_2=\mathit{e}_3, \breve{F}\mathit{e}_3=\mathit{e}_4, \breve{F}\mathit{e}_4=\mathit{e}_2, \breve{F}\mathit{e}_5=\mathit{e}_4, \\
	\breve{F}\mathit{f}_1=\mathit{f}_1, \breve{F}\mathit{f}_2=\mathit{f}_3, \breve{F}\mathit{f}_3=\mathit{f}_4, \breve{F}\mathit{f}_4=\mathit{f}_2, \breve{F}\mathit{f}_5=\mathit{f}_4 
\end{align*}
Let $\vartheta$ be discrete metric on $\breve{E}\cup\breve{P}$ i.e.,
\begin{align*}
	\vartheta(\mathit{e}_\upsilon, \mathit{f}_\varrho)=
\begin{cases}
	1, \quad\text{if}\quad\upsilon\neq\varrho\\
	0, \quad\text{if}\quad\upsilon=\varrho
\end{cases}.
\end{align*}
Consider the mapping $\mathit{q}_0\colon\breve{E}\times\breve{P}\rightarrow[0, \infty)$ defined by
\begin{align*}
	&\mathit{q}_0(\mathit{e}_\upsilon, \mathit{f}_\varrho)=\mathit{q}_0(\mathit{f}_\varrho, \mathit{e}_\upsilon)\\
	&\mathit{q}_0(\mathit{e}_\upsilon, \mathit{e}_\upsilon)=0\\
	&\mathit{q}_0(\mathit{e}_1, \mathit{f}_2)=\mathit{q}_0(\mathit{e}_2, \mathit{f}_3)=4,\\
	&\mathit{q}_0(\mathit{e}_1, \mathit{f}_3)=\mathit{q}_0(\mathit{e}_3, \mathit{f}_4)=3,\\
	&\mathit{q}_0(\mathit{e}_2, \mathit{f}_5)=\mathit{q}_0(\mathit{e}_3, \mathit{f}_5)=5,\\
	&\mathit{q}_0(\mathit{e}_1, \mathit{f}_4)=2,\\
	&\mathit{q}_0(\mathit{e}_1, \mathit{f}_5)=1,\\
	&\mathit{q}_0(\mathit{e}_2, \mathit{f}_4)=7,\\
	&\mathit{q}_0(\mathit{e}_4, \mathit{f}_5)=8.
\end{align*}
We claim that
\begin{align}
	\mathit{q}_0(\breve{F}\mathit{e}, \breve{F}\mathit{f})+\vartheta(\breve{F}\mathit{e}, \breve{F}\mathit{f})\leq\frac{1}{2}\mathit{q}_0(\mathit{e}, \mathit{f})+\vartheta(\mathit{e}, \mathit{f}), \label{eq13}
\end{align}
for every $\mathit{e}\in\breve{E}, \mathit{f}\in\breve{P}$ i.e., $\breve{F}$ is a PC in the sense of definition of PC with $\sigma=1, \mathit{q}_1\equiv1$ and $\pi=\frac{1}{2}$. If $\mathit{e}=\mathit{f}$ (or) $(\mathit{e}, \mathit{f})=(\mathit{e}_1, \mathit{e}_5)$ then equation \eqref{eq13} is obvious. We have to show that equation \eqref{eq13} holds for all $\mathit{e}_\upsilon\in\breve{E}, \mathit{f}_\upsilon\in\breve{P}$ with $1\leq\upsilon<\varrho\leq5$ and $(\upsilon, \varrho)\neq(1,5)$. Then all condition as follows in table below:


\begin{table}[h]
    \centering
    \renewcommand{\arraystretch}{1.2}
    \begin{tabular}{|c|c|c|}
        \hline
        $(\upsilon, \varrho)$ & 
        $\mathit{q}_0(\breve{F}e_\upsilon, \breve{F}f_\varrho) + \vartheta(\breve{F}e_\upsilon, \breve{F}f_\varrho)$ & 
        $\mathit{q}_0(e_\upsilon, f_\varrho) + \vartheta(e_\upsilon, f_\varrho)$ \\
        \hline
        (1,2) & 4 & 5 \\ \hline
        (1,3) & 3 & 4 \\ \hline
        (1,4) & 5 & 3 \\ \hline
        (1,5) & 3 & 2 \\ \hline
        (2,3) & 4 & 5 \\ \hline
        (2,4) & 5 & 8 \\ \hline
        (2,5) & 4 & 6 \\ \hline
        (3,4) & 8 & 4 \\ \hline
        (3,5) & 1 & 6 \\ \hline
        (4,5) & 8 & 9 \\ \hline
        \hline
    \end{tabular}
    \label{tab:bipolar-values}
\end{table}

of Theorem \ref{thm3.1} are fulfilled (\eqref{thm3c2} is satisfied with $\breve{Q}=1$).\\

Hence $\breve{F}$ admits a UFP.
\end{example}

\begin{defn}
	Let $(\breve{E}, \breve{P}, \vartheta)$ be a bipolar $\mathcal{MS}$. The mapping $\breve{F}\colon\breve{E}\cup\breve{P}\rightarrow\breve{E}\cup\breve{P}$ is reffered as picard-continuous, if for all $\mathit{g}\in\breve{E}, \mathit{h}\in\breve{P}$, we have
\begin{align*}
	\lim_{\kappa\rightarrow\infty}\vartheta(\breve{F}^\kappa\mathit{g}, \mathit{h})=0\rightarrow\lim_{\kappa\rightarrow\infty}\vartheta(\breve{F}(\breve{F}^\kappa\mathit{g}), \breve{F}\mathit{h})=0
\end{align*}
where $\breve{F}^0\mathit{g}=\mathit{g}$ and $\breve{F}^{\kappa+1}\mathit{g}=\breve{F}(\breve{F}^\kappa\mathit{g}), \forall\kappa\geq0$.
\end{defn}
\begin{remark}
	Remark that, if $\breve{F}\colon\breve{E}\cup\breve{P}\rightarrow\breve{E}\cup\breve{P}$ is continuous, then $\breve{F}$ is picard-continuous, the converse part is not true. 
\end{remark}
\begin{example}
	Let $\breve{E}=[\mathit{q}, \mathit{w}], \breve{P}=[\mu, \vartheta]$, where $\mathit{q}, \mathit{w}, \mu, \vartheta\in\mathbb{R}$ and $\mathit{q}<\mathit{w}, \mu<\vartheta$. The mapping $\breve{F}\colon\breve{E}\cup\breve{P}\rightarrow\breve{E}\cup\breve{P}$ defined by
\begin{align*}
	\breve{F}\mathit{e}&=
\begin{cases}
	\mathit{q}, \quad\text{if}\quad\mathit{q}\leq\mathit{e}<\mathit{w},\\
	\frac{\mathit{q}+\mathit{w}}{3}, \quad\text{if}\quad\mathit{e}=\mathit{w}
\end{cases}\\
\breve{F}\mathit{f}&=
\begin{cases}
	\mu, \quad\text{if}\quad\mu\leq\mu<\vartheta,\\
	\frac{\mu+\vartheta}{3}, \quad\text{if}\quad\mathit{f}=\vartheta.
\end{cases}.
\end{align*}
Let $\vartheta$ be the standard metric on $\breve{E}$, i.e., $\vartheta(\mathit{e}, \mathit{f})=|\mathit{e}-\mathit{f}|$ forall $\mathit{e}, \mathit{f}\in\breve{E}\cup\breve{P}$. Although the mapping $\breve{F}$ is not continuous at 
$\mathit{w}$ and $\vartheta$, it is Picard-continuous as per the preceding definition.
\begin{align*}
	\breve{F}^\kappa\mathit{g}=\mathit{q}, \quad\forall\kappa\geq3.
\end{align*}
So, if for some $\mathit{g}\in\breve{E}, \mathit{h}\in\breve{P}$, we have
\begin{align*}
	\lim_{\kappa\rightarrow\infty}\vartheta(\breve{F}^\kappa\mathit{g}, \mathit{h})=0,
\end{align*}
then $\mathit{h}=\mathit{q}$.
\begin{align*}
	\lim_{\kappa\rightarrow\infty}\vartheta(\breve{F}(\breve{F}^\kappa\mathit{g}), \mathit{h})=\lim_{\kappa\rightarrow\infty}\vartheta(\breve{F}\mathit{q}, \mathit{h})=\vartheta(\breve{F}\mathit{q}, \mathit{h})=0.
\end{align*}
Hence $\breve{F}$ is picard-continuous.
\end{example}
\begin{thm}
	Let $(\breve{E}, \breve{P}, \vartheta)$ be a complete bipolar $\mathcal{MS}$ and the mapping $\breve{F}\colon\breve{E}\cup\breve{P}\rightarrow\breve{E}\cup\breve{P}$ be a PC such that
\begin{enumerate}[label=(\roman*)]
	\item $\breve{F}$ is picard-continuous,
	\item We can find that $\varrho\in\{1,\cdots,\sigma\}$ and $\breve{Q}_\varrho>0$ implies that
\begin{align*}
	\mathit{q}_\varrho(\mathit{e}, \mathit{f})\geq\breve{Q}_\varrho, \quad\mathit{e}\in\breve{E}, \mathit{f}\in\breve{P}.
\end{align*}
\end{enumerate}
Then $\breve{F}$ admits a UFP.
\begin{proof}
	Let $\breve{F}$ be a non-empty set. Assume $\mathit{g}_0\in\breve{E}, \mathit{h}_0\in\breve{P}$ and $\{\mathit{g}_\kappa\}\subset\breve{E}$ and $\{\mathit{h}_\kappa\}\subset\breve{P}$ be the picard bisequence defined by
\begin{align*}
	\mathit{g}_{\kappa+1}=\breve{F}\mathit{g}_\kappa \quad\text{and}\quad\mathit{h}_{\kappa+1}=\breve{F}\mathit{h}_\kappa, \quad\forall\kappa\geq0,
\end{align*}
i.e.,
\begin{align*}
	\mathit{g}_\kappa=\breve{F}^\kappa\mathit{g}_0\quad\text{and}\quad\mathit{h}_\kappa=\breve{F}^\kappa\mathit{h}_0,\quad\forall\kappa\geq0.
\end{align*}
By using the proof of Theorem \ref{thm3.1}, w.k.t $\{\mathit{g}_\kappa\}$ and $\{\mathit{h}_\kappa\}$ is a Cauchy bisequence, by completeness of $(\breve{E}, \breve{P}, \vartheta)$ that we can find that $\mathit{g}\in\breve{E}$ and $\mathit{h}\in\breve{P}$ implies that
\begin{align*}
	\lim_{\kappa\rightarrow\infty}\vartheta(\breve{F}^\kappa\mathit{g}_0, \mathit{g})=0 \quad\text{and}\quad\lim_{\kappa\rightarrow\infty}\vartheta(\breve{F}^\kappa\mathit{h}_0, \mathit{h})=0.
\end{align*}
Using the picard continuity of $\breve{F}$, it hold that
\begin{align*}
	\lim_{\kappa\rightarrow\infty}\vartheta(\breve{F}^{\kappa+1}\mathit{g}_0, \breve{F}\mathit{g})=\lim_{\kappa\rightarrow\infty}\vartheta(\breve{F}(\breve{F}^\kappa\mathit{g}_0), \breve{F}\mathit{g})=0,
\end{align*}
and
\begin{align*}
	\lim_{\kappa\rightarrow\infty}\vartheta(\breve{F}^{\kappa+1}\mathit{h}_0, \breve{F}\mathit{h})=\lim_{\kappa\rightarrow\infty}\vartheta(\breve{F}(\breve{F}^\kappa\mathit{h}_0), \breve{F}\mathit{h})=0,
\end{align*}
using the uniqueness of limit, then $\mathit{g}$ is a FP of $\breve{F}$.
\end{proof}
\end{thm}
\section{The class of almost polynomial contraction}
\begin{defn}
	Let $(\breve{E}, \breve{P}, \vartheta)$ be a bipolar $\mathcal{MS}$ and the mapping $\breve{F}\colon\breve{E}\cup\breve{P}\rightarrow\breve{E}\cup\breve{P}$. We say that $\breve{F}$ is an almost PC, if we can find that $\pi\in(0, 1)$,  $\sigma\geq1$, a sequence $\{\breve{H}_\upsilon\}_{\upsilon=0}^{\sigma}\subset(0, \infty)$ and a mapping $\mathit{q}_\upsilon\colon\breve{E}\times\breve{P}\rightarrow[0, \infty)$, $\upsilon=0,\cdots,\sigma$  implies that
\begin{align*}
	\sum_{\upsilon=0}^{\sigma}\mathit{q}_\upsilon(\breve{F}\mathit{e}, \breve{F}\mathit{f})\vartheta^\upsilon(\breve{F}\mathit{e}, \breve{F}\mathit{f})\leq\pi\sum_{\upsilon=0}^{\sigma}\mathit{q}_\upsilon(\mathit{e}, \mathit{f})[\vartheta^\upsilon(\mathit{e}, \mathit{f})+\breve{H}_\upsilon\vartheta^\upsilon(\mathit{f}, \breve{F}\mathit{e})], \quad\mathit{e}\in\breve{E}, \mathit{f}\in\breve{P}.
\end{align*}
\end{defn}
\begin{defn}
	Let $(\breve{E}, \breve{P}, \vartheta)$ be a bipolar $\mathcal{MS}$ and the mapping be $\breve{F}\colon\breve{E}\cup\breve{P}\rightarrow\breve{E}\cup\breve{P}$. We reffered as $\breve{F}$ is a weakly picard operator, if
\begin{enumerate}[label=(\roman*)]
	\item All FPs of $\breve{F}$ is non-empty,
	\item Forall $\mathit{g}_0\in\breve{E}, \mathit{h}_0\in\breve{P}$, the picard sequence $\{\breve{F}^\kappa\mathit{g}_0\}$ and $\{\breve{F}^\kappa\mathit{h}_0\}$ is biconvergent and its limit belongs to the set of FPs of $\breve{F}$.
\end{enumerate}
\end{defn}
\begin{thm} \label{thm4.1}
	Let $(\breve{E}, \breve{P}, \vartheta)$ be a complete bipolar $\mathcal{MS}$ and the covariant mapping $\breve{F}\colon\breve{E}\cup\breve{P}\rightrightarrows\breve{E}\cup\breve{P}$ be an atmost PC such that
\begin{enumerate}[label=(\roman*)]
	\item $\breve{F}$ is picard-continuous,
	\item we can find that $\varrho\in\{1,\cdots,\sigma\}$ and $\breve{Q}_\varrho>0$ implies that \label{thm4.1c2}
\begin{align*}
	\mathit{q}_\varrho(\mathit{e}, \mathit{f})\geq\breve{Q}_\varrho, \quad\mathit{e}\in\breve{E}, \mathit{f}\in\breve{P}.
\end{align*}
\end{enumerate}
Then $\breve{F}$ is a weakly picard operator.
\begin{proof}
	Let $\mathit{g}_0\in\breve{E}$ and $\mathit{h}_0\in\breve{P}$ be fixed and $\{\mathit{g}_\kappa\}\subset\breve{E}$ and $\{\mathit{h}_\kappa\}\subset\breve{P}$ be the picard bisequence defined by
\begin{align*}
	\mathit{g}_{\kappa+1}=\breve{F}\mathit{g}_\kappa\quad\text{and}\quad\mathit{h}_{\kappa+1}=\breve{F}\mathit{h}_\kappa, \quad\forall\kappa\geq0.
\end{align*}
By the definition of PC with $(\mathit{e}, \mathit{f})=(\mathit{g}_0, \mathit{h}_0)$
\begin{align*}
	\sum_{\upsilon=0}^{\sigma}\mathit{q}_\upsilon(\breve{F}\mathit{g}_0, \breve{F}\mathit{h}_0)\vartheta^\upsilon(\breve{F}\mathit{g}_0, \breve{F}\mathit{h}_0)\leq\pi\sum_{\upsilon=0}^{\sigma}\mathit{q}_\upsilon(\mathit{g}_0, \mathit{h}_0)[\vartheta^\upsilon(\mathit{g}_0, \mathit{h}_0)+\breve{H}_\upsilon\vartheta^\upsilon(\mathit{h}_0, \breve{F}\mathit{g}_0)],
\end{align*}
i.e.,
\begin{align*}
	\sum_{\upsilon=0}^{\sigma}\mathit{q}_\upsilon(\mathit{g}_1, \mathit{h}_1)\vartheta^\upsilon(\mathit{g}_1, \mathit{h}_1)\leq\pi\sum_{\upsilon=0}^{\sigma}\mathit{q}_\upsilon(\mathit{g}_0, \mathit{h}_0)\vartheta^\upsilon(\mathit{g}_0, \mathit{h}_0).
\end{align*}
Again by definition of PC with $(\mathit{e}, \mathit{f})=(\mathit{g}_1, \mathit{h}_1)$, we get
\begin{align*}
	\sum_{\upsilon=0}^{\sigma}\mathit{q}_\upsilon(\breve{F}\mathit{g}_1, \breve{F}\mathit{h}_1)\vartheta^\upsilon(\breve{F}\mathit{g}_1, \breve{F}\mathit{h}_1)\leq\pi\sum_{\upsilon=0}^{\sigma}\mathit{q}_\upsilon(\mathit{g}_1, \mathit{h}_1)[\vartheta^\upsilon(\mathit{g}_1, \mathit{h}_1)+\breve{H}_\upsilon\vartheta^\upsilon(\mathit{h}_1, \breve{F}\mathit{g}_1)],
\end{align*}
i.e.,
\begin{align*}
	\sum_{\upsilon=0}^{\sigma}\mathit{q}_\upsilon(\mathit{g}_2, \mathit{h}_2)\vartheta^\upsilon(\mathit{g}_2, \mathit{h}_2)\leq\pi\sum_{\upsilon=0}^{\sigma}\mathit{q}_\upsilon(\mathit{g}_1, \mathit{h}_1)\vartheta^\upsilon(\mathit{g}_1, \mathit{h}_1),
\end{align*}
which implies
\begin{align*}
	\sum_{\upsilon=0}^{\sigma}\mathit{q}_\upsilon(\mathit{g}_2, \mathit{h}_2)\vartheta^\upsilon(\mathit{g}_2, \mathit{h}_2)\leq\pi\sum_{\upsilon=0}^{\sigma}\mathit{q}_\upsilon(\mathit{g}_0, \mathit{h}_0)\vartheta^\upsilon(\mathit{g}_0, \mathit{h}_0),
\end{align*}
continuing by induction process, we get
\begin{align*}
	\sum_{\upsilon=0}^{\sigma}\mathit{q}_\upsilon(\mathit{g}_\kappa, \mathit{h}_\kappa)\vartheta^\upsilon(\mathit{g}_\kappa, \mathit{h}_\kappa)\leq\pi^\kappa\sum_{\upsilon=0}^{\sigma}\mathit{q}_\upsilon(\mathit{g}_0, \mathit{h}_0)\vartheta^\upsilon(\mathit{g}_0, \mathit{h}_0), \quad\kappa\geq0,
\end{align*}
which implies from \ref{thm4.1c2} that
\begin{align*}
	\vartheta^\varrho(\mathit{g}_\kappa, \mathit{h}_\kappa)\leq\pi^\kappa\mathcal{M}, \quad\kappa\geq0,
\end{align*}
where $\mathcal{M}$ is given in \eqref{eq4}. By the proof of Theorem \ref{thm3.1}. We obtain that $\{\mathit{g}_\kappa\}$ and $\{\mathit{h}_\kappa\}$ is a cauchy bisequence, by the completeness of $(\breve{E}, \breve{P}, \vartheta)$ the existence of $\mathit{g}\in\breve{E}$ implies that
\begin{align*}
	\lim_{\kappa\rightarrow\infty}\vartheta(\mathit{g}_\kappa, \mathit{g})=0.
\end{align*}
In conclusion, considering that $\breve{F}$ is picard-continuous, we obtain
\begin{align*}
	\lim_{\kappa\rightarrow\infty}\vartheta(\mathit{g}_{\kappa+1}, \breve{F}\mathit{g})=0,
\end{align*}
by the uniqueness of limit that $\mathit{g}=\breve{F}\mathit{g}$. This completes the proof.
\end{proof}
\end{thm}
\begin{prop}
	Let $(\breve{E}, \breve{P}, \vartheta)$ be a bipolar $\mathcal{MS}$ and the mapping $\breve{F}\colon\breve{E}\cup\breve{P}\rightarrow\breve{E}\cup\breve{P}$. Assume that there exist $\pi\in(0, 1)$,  $\sigma\geq1$ and two finite sequence $\{\mathit{q}_\upsilon\}_{\upsilon=1}^{\sigma}, \{\breve{H}_\upsilon\}_{\upsilon=1}^{\sigma}\subset(0, \infty)$ such that
\begin{align}
	\sum_{\upsilon=1}^{\sigma}\mathit{q}_\upsilon\vartheta^\upsilon(\breve{F}\mathit{e}, \breve{F}\mathit{f})\leq\pi\sum_{\upsilon=1}^{\sigma}\mathit{q}_\upsilon[\vartheta^\upsilon(\mathit{e}, \mathit{f})+\breve{H}_\upsilon\vartheta^\upsilon(\mathit{f}, \breve{F}\mathit{e})], \label{eq14}
\end{align}
for every $\mathit{e}\in\breve{E}, \mathit{f}\in\breve{P}$. Then $\breve{F}$ is picard-continuous.
\begin{proof}
	Let $\mathit{g}\in\breve{E}, \mathit{f}\in\breve{P}$ be such that
\begin{align*}
	\lim_{\kappa\rightarrow\infty}\vartheta(\breve{F}^\kappa\mathit{g}, \mathit{h})=0.
\end{align*}
By equation \eqref{eq14} with $(\mathit{e}, \mathit{f})=(\breve{F}^\kappa\mathit{g}, \mathit{h})$, we get
\begin{align*}
	\sum_{\upsilon=1}^{\sigma}\mathit{q}_\upsilon\vartheta^\upsilon(\breve{F}(\breve{F}^\kappa\mathit{g}), \breve{F}\mathit{h})\leq\pi\sum_{\upsilon=1}^{\sigma}\mathit{q}_\upsilon[\vartheta^\upsilon(\breve{F}^\kappa\mathit{g}, \mathit{h})+\breve{H}_\upsilon\vartheta^\upsilon(\mathit{h}, \breve{F}(\breve{F}^\kappa\mathit{g}))],
\end{align*}
i.e.,
\begin{align*}
	\sum_{\upsilon=1}^{\sigma}\mathit{q}_\upsilon\vartheta^\upsilon(\breve{F}^{\kappa+1}\mathit{g}, \breve{F}\mathit{h})\leq\pi\sum_{\upsilon=1}^{\sigma}\mathit{q}_\upsilon[\vartheta^\upsilon(\breve{F}^\kappa\mathit{g}, \mathit{h})+\breve{H}_\upsilon\vartheta^\upsilon(\mathit{h}, \breve{F}^{\kappa+1}\mathit{g})],
\end{align*}
which implies that
\begin{align*}
	\vartheta(\breve{F}(\breve{F}^\kappa\mathit{g}), \breve{F}\mathit{h})\leq\frac{\pi}{\mathit{q}_1}\sum_{\upsilon=1}^{\sigma}\mathit{q}_\upsilon[\vartheta^\upsilon(\breve{F}^\kappa\mathit{g}, \mathit{h})+\breve{H}_\upsilon\vartheta^\upsilon(\mathit{h}, \breve{F}^{\kappa+1}\mathit{g})].
\end{align*}
Then, passing to the limit as $\kappa\rightarrow\infty$ in the above inequality and by equation \eqref{eq9}, we get
\begin{align*}
	\lim_{\kappa\rightarrow\infty}\vartheta(\breve{F}(\breve{F}^\kappa\mathit{g}), \breve{F}\mathit{h})=0.
\end{align*}
Hence $\breve{F}$ is picard-continuous.
\end{proof}
\end{prop}
The following result is derived from the above theorem and proposition.
\begin{cor}
	Let $(\breve{E}, \breve{P}, \vartheta)$ be a complete bipolar $\mathcal{MS}$ and the mapping $\breve{F}\colon\breve{E}\cup\breve{P}\rightarrow\breve{E}\cup\breve{P}$. Assume that there exist $\pi\in(0, 1)$, $\sigma\geq1$, and two finite sequence $\{\mathit{q}_\upsilon\}_{\upsilon=1}^{\sigma},\{\breve{H}_\upsilon\}_{\upsilon=1}^{\sigma}\subset(0, \infty)$ such that equation \eqref{eq9} holds for every $\mathit{e}\in\breve{E}, \mathit{f}\in\breve{P}$. Then, $\breve{F}$ is a weakly picard continuous.
\begin{proof}
	By equation \eqref{eq9} is a special case by defintion of almost PC with $\mathit{q}_0\equiv0$ and $\mathit{q}_\upsilon$ is constant, $\forall\upsilon\in\{1,\cdots,\sigma\}$. Then by above proposition and theorem applies.
\end{proof}
\end{cor}
\begin{remark}
	Taking $\sigma=1, \mathit{q}_1=1$ and $\breve{H}_1=\frac{\rho}{\pi}$, where $\rho>0$, equation \eqref{eq9} reduces to PC. Then, by above corollary, we recover Berinde's FPT.
\end{remark}
\begin{example}
	Let $\breve{E}=\{\mathit{e}_1, \mathit{e}_2\}, \breve{P}=\{\mathit{e}_1, \mathit{f}_2\}$ and the mapping be $\breve{F}\colon\breve{E}\cup\breve{P}\rightarrow\breve{E}\cup\breve{P}$ defined by
\begin{align*}
	\breve{F}\mathit{e}_1=\mathit{e}_1, \breve{F}\mathit{e}_2=\mathit{e}_1, \breve{F}\mathit{f}_1=\mathit{e}_1, \breve{F}\mathit{f}_2=\mathit{e}_1.
\end{align*}

Define the two metrics as
\begin{align*}
	\vartheta_1(\mathit{e}_\upsilon, \mathit{f}_\varrho)&=
\begin{cases}
	1 &\mathit{e}_\upsilon\neq\mathit{f}_\varrho\\
	0 &\mathit{e}_\upsilon=\mathit{f}_\varrho
\end{cases}\\
\vartheta_2(\mathit{e}_\upsilon, \mathit{f}_\varrho)&=
\begin{cases}
	2 &\mathit{e}_\upsilon\in\breve{E}, \mathit{f}_\varrho\in\breve{P}\\
	1 &\mathit{e}_\upsilon\neq\mathit{f}_\varrho\\
	0 &\mathit{e}_\upsilon=\mathit{f}_\varrho
\end{cases}
\end{align*}

\begin{table}[h]
	\centering
    \renewcommand{\arraystretch}{1.2}
	\begin{tabular}{|c|c|c|c|c|}
		\hline
		$(\mathit{e}_\upsilon, \mathit{f}_\varrho)$ & $\vartheta_1(\mathit{e}_\upsilon, \mathit{f}_\varrho)$ & $\vartheta_2(\mathit{e}_\upsilon, \mathit{f}_\varrho)$ &$\vartheta_1(\breve{F}\mathit{e}_\upsilon, \breve{F}\mathit{f}_\varrho)$ &$\vartheta_2(\breve{F}\mathit{e}_\upsilon, \breve{F}\mathit{f}_\varrho)$  \\ \hline
		$(\mathit{e}_1, \mathit{e}_2)$ &1 &1 &0&0	\\ \hline
		$(\mathit{e}_1, \mathit{f}_1)$ &1 &2 &0&0	\\ \hline
		$(\mathit{e}_2, \mathit{f}_2)$ &1 &2 &0&0	\\ \hline
		$(\mathit{f}_1, \mathit{f}_2)$ &1 &1 &0&0	\\ \hline		
	\end{tabular}
\end{table}

From the table:

\begin{align*}
	\vartheta_1(\breve{F}\mathit{e}_\upsilon, \breve{F}\mathit{f}_\varrho)\leq\vartheta_1(\mathit{e}_\upsilon, \mathit{e}_\varrho), \vartheta_2(\breve{F}\mathit{e}_\upsilon, \breve{F}\mathit{f}_\varrho)\leq\vartheta_2(\mathit{e}_\upsilon, \mathit{e}_\varrho) \quad\mathit{e}\in\breve{E}, \mathit{f}\in\breve{P}.
\end{align*}

Hence the mapping $\breve{F}$ is non-expansive with bipolar $\mathcal{MS}$.
\end{example}
\section{Application to Fractional Calculus}
Let $(\breve{V}[0, 1])$ be the set of all continuous function on $[0, 1]$ and the mapping be $\vartheta\colon\breve{V}([0, 1])\times\breve{V}([0, 1])\rightarrow\mathbb{R}$ defined by $\vartheta(\mathit{e}, \mathit{f})=\|\mathit{e}-\mathit{f}\|_\infty=\sup_{\rho\in[0, 1]}|\mathit{e}(\rho)-\mathit{f}(\rho)|$. Define $\mathit{e}\in\breve{E}$ and $\mathit{f}\in\breve{P}$.\\
The Caputo fractional derivative of order $\beta$ applied to a continuous function $\mathfrak{h}\colon[0, +\infty)\rightarrow\mathbb{R}$ is defined as
\begin{align*}
	^{\breve{V}}\breve{O}^\beta(\mathfrak{h}(\rho))=\frac{1}{\Gamma(\varpi-\beta)}\int_{0}^{\rho}(\rho-\eta)^{\varpi-\beta-1}\mathfrak{g}^{(\varpi)}(\eta)\mathit{d}\eta, \quad(\varpi-1<\beta<\kappa, \varpi=[\beta]+1),
\end{align*}
where $\Gamma$ is the gamma function and $[\beta]$ denotes the integer part of a real number. Consider $\breve{E}=(\breve{V}[0,1], [0,\infty)=\{\vartheta:[0,1]\rightarrow[0, \infty)$ be a continuous\} and $\breve{P}=(\breve{V}[0,1], [0,\infty)=\{\vartheta:[0,1]\rightarrow[-\infty,0)$ be a continuous\}. \\
Let us see, the existence of the solution of non-linear fractional differential equation
\begin{align}
	^{\breve{V}}\breve{O}^\beta(\mathit{g}(\rho))+\omega(\rho, \mathit{g}(\rho))=0 \quad(0\leq\rho\leq1, \beta<1), \label{1}
\end{align}
with $\mathit{g}(0)=\mathit{g}(1)=0$ and $\omega\colon[0, 1]\times\mathbb{R}\rightarrow\mathbb{R}$ is a continuous function, and the Green’s function corresponding to problem \eqref{1} is defined as,
\begin{align*}
	\breve{A}(\rho, \eta)=
	\begin{cases}
		\rho(1-\eta)^{\mathit{q}-1}-(\rho-\eta)^{\mathit{q}-1}, &\text{if}\, \, 0\leq\rho\leq\eta\leq1, \\
		\frac{\rho(1-\eta)^{\mathit{q}-1}}{\Gamma(\mathit{q})}, &\text{if}\, \, 0\leq\eta\leq\rho\leq1.
	\end{cases}
\end{align*}
Suppose the following conditions are satisfied:
\begin{enumerate}[label=(\arabic*)]
	\item There exist $\sigma\in(0,1)$ such that $|\omega(\rho, \mathit{e})-\omega(\rho, \mathit{f})|\leq\sigma|\mathit{e}-\mathit{f}|$ for each $\rho\in[0, 1]$ and $\mathit{e}, \mathit{f}\in\mathbb{R}$, \label{con-1}
	\item $\sup_{\rho\in[0, 1]}\bigg(\int_{0}^{1}\breve{A}(\rho, \eta)\mathit{d}\eta\bigg)\leq1$. \label{con-2}
\end{enumerate}
Next, we establish the existence of a solution to the fractional differential equation \eqref{1}.
\begin{thm}
	Assuming the condition holds \eqref{con-1}-\eqref{con-2}, from \eqref{1} has a unique solution.
	\begin{proof}
		Define the mapping $\breve{F}\colon\breve{V}[0, 1]\rightarrow\breve{V}[0, 1]$ is defined by
		\begin{align*}
			\breve{F}(\mathit{e}(\rho))=\int_{0}^{1}\breve{A}(\rho, \eta)\omega(\eta, \mathit{e}(\eta))\mathit{d}\eta.
		\end{align*}
		Under our assumption, $\mathit{e}$ solves \eqref{1} precisely when $\mathit{e}\in\breve{E}$ solves the integral equation.
		\begin{align*}
			\mathit{g}(\rho)=\int_{0}^{1}\breve{A}(\rho, \eta)\omega(\eta, \mathit{g}(\eta))\mathit{d}\eta, \, \, \forall\rho\in[0, 1].
		\end{align*}
		Consider,
		\begin{align*}
			|\breve{F}\mathit{e}(\rho)-\breve{F}\mathit{f}(\rho)|&=\bigg|\int_{0}^{1}\breve{A}(\rho, \eta)\omega(\eta, \mathit{e}(\eta))\mathit{d}\eta-\int_{0}^{1}\breve{A}(\rho, \eta)\omega(\eta, \mathit{f}(\eta))\mathit{d}\eta\bigg|\\
			&\leq\int_{0}^{1}|\breve{A}(\rho, \eta)(\omega(\eta, \mathit{e}(\eta))-\omega(\eta, \mathit{f}(\eta)))\mathit{d}\eta|\\
			&\leq\int_{0}^{1}\breve{A}(\rho, \eta) |(\omega(\eta, \mathit{e}(\eta))-\omega(\eta, \mathit{f}(\eta)))|\mathit{d}\eta\\
			&\leq\int_{0}^{1}\breve{A}(\rho, \eta)\sigma |(\mathit{e}(\eta)- \mathit{f}(\eta))|\mathit{d}\eta.
		\end{align*}
		Taking supremum on both sides, we get 
		\begin{align*}
			\vartheta(\breve{F}\mathit{e}, \breve{F}\mathit{f})\leq\sigma\vartheta(\mathit{e}, \mathit{f}).
		\end{align*}
		Therefore, all the hypothesis of theorem \ref{thm3.1} and \ref{thm4.1}  are satisfied. Hence, $\breve{F}$ has a unique solution in $\breve{V}[0, 1]$. 
	\end{proof}
\end{thm}

\section{Conclusion}
In this work, we have explored FP results within the framework of bipolar $\mathcal{MS}$s by introducing and analyzing polynomial-type contractions and almost polynomial contractions. The generalization of classical contraction principles under this extended metric structure offers a richer context for addressing FP problems. The results obtained not only unify and extend several existing theorems in the literature but also establish a novel approach applicable to fractional differential equations. The applicability of the developed FPT to problems in fractional calculus highlights the potential of bipolar $\mathcal{MS}$s as a powerful tool in nonlinear analysis. Future research may focus on extending these results to multivalued mappings, dynamic systems, and other abstract spaces, thereby broadening the theoretical and practical impact of the developed concepts.

	\bibliographystyle{amsplain}

\begin{thebibliography}{20}
	\bibitem{01} P. Agarwal, M. Jleli, B. Samet, Fixed Point Theory in Metric Spaces, Recent Advances and Applications, Springer, Singapore (2018).
	
	\bibitem{02} S. Banach, Sur les op´erations dans les ensembles abstraits et leur applications aux
	´equations int´egrales, Fund Math. 3 (1922) 133–181.
	
	\bibitem{03} V. Berinde, Approximating fixed points of weak contractions using the Picard iteration,	Nonlinear Analysis Forum. 9(1) (2004) 43–53.
	
	\bibitem{04} V. Berinde, General constructive fixed point theorems for \'{C}iri\'{c}-type almost contractions
	in metric spaces, Carpathian J. Math. 24 (2008) 10–19.
	
	\bibitem{05} V. Berinde, Approximating fixed points of enriched nonexpansive mappings in Banach
	spaces by using a retraction-displacement condition, Carpathian J. Math. 36(1) (2020)
	27–34.
	
	\bibitem{06} V. Berinde, M. Pˇacurar, Fixed point theorems for nonself single-valued almost contractions, Fixed Point Theory. 14 (2013) 301–311.
	
	\bibitem{07} V. Berinde, M. P\u{a}curar, Approximating fixed points of enriched contractions in Banach
	spaces, J. Fixed Point Theory Appl. 22 (2020) 1–10.
	
	\bibitem{08} DW. Boyd, JSW. Wong, On nonlinear contractions, Proc Am Math Soc. 20 (1969)
	458–464.

    \bibitem{g1} Y.U.Gaba, M.Aphane, V.Sihag, On Two Banach‑Type Fixed Points in Bipolar Metric Spaces, Abstract and Applied Analysis, vol.2021 (2021) 1–10.

    \bibitem{g2} Y.U.Gaba, M.Aphane, H.Aydi, $(\alpha,BK)$-Contractions in Bipolar Metric Spaces, Journal of Mathematics, vol.2021 (2021) 1–6.


	
	\bibitem{09} A. Branciari, A fixed point theorem of Banach-Caccioppoli type on a class of generalized	metric spaces, Publ. Math. Debr. 57 (2000) 31–37.
	
	\bibitem{10} S.K. Chatterjea, Fixed-point theorems, C.R. Acad. Bulgare Sci. 25 (1972) 727–730.
	
	\bibitem{11} Lj. B. \'{C}iri\'{c}, A generalization of Banach’s contraction principle, Proc. Amer. Math. Soc.	45 (2) (1974) 267–273.
	
	\bibitem{12} S. Czerwik, Contraction mappings in b-metric spaces, Acta Math. Inform., Univ. Ostrav.
	1(1) (1993) 5–11.
	
	\bibitem{13} B.C. Dhage, Generalized metric space and mapping with fixed point, Bulletin of the
	Calcutta Mathematical Society. 84 (1992) 329–336.
	
	\bibitem{14} E. Petrov, Fixed point theorem for mappings contracting perimeters of triangles. J.
	Fixed Point Theory Appl. 25 (2023) 74.
	
	\bibitem{15} A. Petru¸sel, I.A. Rus, Fixed point theory in terms of a metric and of an order relation,Fixed Point Theory. 20(2) (2019) 601–622.
	
	\bibitem{16} O. Popescu, C. P\u{a}curar, Fixed point theorem for generalized Chatterjea type mappings, arXiv:2404.00782.
	
	\bibitem{17} E. Rakotch, A note on contractive mappings, Proc Am Math Soc. 13 (1962) 459–465.
	
	\bibitem{18} S. Reich, Fixed points of contractive functions, Boll Unione Mat Ital. 5 (1972) 26–42.
	
	\bibitem{19} I.A. Rus, Weakly Picard mappings, Comment. Mat. Univ. Carolinae. 34(4) (1993) 769–
	773.
	
	\bibitem{20} I.A. Rus, Generalized Contractions and Applications, Cluj University Press, Cluj-
	Napoca, 2001.
	
	\bibitem{21} I.A. Rus, Picard operators and applications, Scientiae Math. Japon. 58(1) (2003) 191–219.
	
	\bibitem{22} I.A. Rus, A. Petru\c{s}el, M.A. \c{S}erban, Weakly Picard operators: equivalent definitions, applications and open problems, Fixed Point Theory. 7 (2006) 3–22.
	
	\bibitem{23} Q. Zhang, Y. Song, Fixed point theory for generalized $(\psi-\varphi)$-weak contractions, Appl.Math. Lett. 22 (2009) 75–78.
	
	\bibitem{24} G. Mani, R. Ramaswamy, A. J. Gnanaprakasam, V. Stojiljkovi\'{c}, Z.M. Fadail, S. Radenovi\'{c}, Application of fixed point results in the setting of 
	$F$-contraction and simulation function in the setting of bipolar metric space[J]. AIMS Mathematics, 2023, 8(2): 3269-3285. doi: 10.3934/math.2023168.
	
	\bibitem{25} G. Mani,  A.J. Gnanaprakasam,  Z.D. Mitrovi\'{c},  M.F. Bota, Solving an Integral Equation via Fuzzy Triple Controlled Bipolar Metric Spaces. Mathematics 2021, 9, 3181. https://doi.org/10.3390/math9243181
	
	\bibitem{26} G. Mani, G. Janardhanan, S.T.M. Thabet, I. Kedim, M.V. Cortez, Common Fixed Point Techniques in Bipolar Orthogonal Metric Space WithApplications to Economic Problem and Integral Equation, Int. J. Anal. Appl. (2025), https://doi 10.28924/2291-8639-23-2025-52
	
	\bibitem{27} T. Abdeljawad, N. Khan, B. Abdalla,  A. Al-Jaser, M. Alqudah, \& K. Shah, A mathematical analysis of human papilloma virus (HPV) disease with new perspectives of fractional calculus. Alexandria Engineering Journal, (2025), 126, 575–599. https://doi.org/10.1016/j.aej.2025.03.136
	
	\bibitem{28}R. Begum,  S. Ali, T. Abdeljawad, \& K. Shah,  The self-protection behavior changes effect in SEIR-DM of COVID-19 pandemic model with numerical simulation and controlling strategies presentation. Fractals, (2025), 33(02), 2540140. https://doi.org/10.1142/S0218348X25401401
	
	\bibitem{29} A.A. Thirthar,  H. Abboubakar,  A.A. Lamrani, K.S. Nisar, Dynamical behavior of a
	fractional-order epidemic model for investigating two fear effect functions. Results in
	Control and Optimization. 2024, 16, 2666-7207.
	
	\bibitem{30} K. Muthuvel, K. Kaliraj,  S.N. Kottakkaran, V. Vijayakumar, Relative controllability
	for $\psi$-Caputo fractional delay control system. Results in Control and Optimization.
	2024, 16, 2666-7207.
	
	\bibitem{31} K.S. Nisar, A constructive numerical approach to solve the Fractional Modified
	Cassama-Holm equation. Alexandria Engineering. 2024, 106, 19-24.
   
	\end{thebibliography}

\end{document}